\documentclass[12pt]{article}

\usepackage {amsfonts,amsthm}
\usepackage[cp1251]{inputenc}
\usepackage[english]{babel}
\theoremstyle{plain}
\newtheorem{Th}{Theorem}

\newtheorem{Def}{Definition}
\def\o{\omega}
\def\bs{~\hfill\rule{7pt}{7pt}}
\def\R{{\mathbb R}}
\def\Tra{{\rm Tra}}
\def\card{{\rm card}}
\def\supp{{\rm supp}}
\def\N{{\mathbb N}}
\def\Z{{\mathbb Z}}
\def\a{\alpha}
\def\b{\beta}
\def\g{\gamma}
\def\l{\lambda}

\def\p{\varphi}

\begin{document}

\title{
    \bf
A uniformly spread measure criterion}

\author{A.\,Dudko, S.\,Favorov}

\date{}
\maketitle

\begin{abstract}
We prove that if all shifts of a measure in the Euclidean space are close in a sense
to each  other, then this measure is close to the Lebesgue one.
\end{abstract}

Let $(a_n)_{n\in\N}$ be a discrete sequence in $\R^d$, i.e., a map $\N\to\R^d$ such
that its image $\{a_n\}_{n\in\N}$ has no limit point in $\R^d$ and each
$x\in\{a_n\}_{n\in\N}$ has at most a finite multiplicity. Following Laczkovich
\cite{L1}, \cite{L2}, we say that the sequence is uniformly spread over $\R^d$, if
there is $\a>0$  such that
\begin{equation}\label{d0}
\inf_\psi\sup_n|a_n-\a\psi(n)|<\infty,
\end{equation}
where the infimum  is taken over all bijections $\psi:\,\N\to\Z^d$.

Using the idea of the mass transfer (see, for example, \cite{KA}), M.~Sodin,
B.~Tsirelson  extended the above definition to  measures in $\R^d$. In \cite{ST} they
introduce {\it a transportation distance} between  arbitrary locally finite positive
measures $\nu_1$ and  $\nu_2$
$$
\Tra(\nu_1,\nu_2)=\inf_\g\sup\{|x-y|:\,x,y\in\overline{\supp\g}\}.
 $$
Here the infimum  is taken over all {\it transportation measures} $\g$ between
measures $\nu_1$ and $\nu_2$, where the latter means
 \begin{equation}\label{t1}
\int\int_{\R^d\times\R^d}\p(x)~d\g(x,y)=\int_{\R^d}\p(x)~d\nu_1(x),
 \end{equation}
\begin{equation}\label{t2}
\int\int_{\R^d\times\R^d}\p(y)~d\g(x,y)=\int_{\R^d}\p(y)~d\nu_2(y)
 \end{equation}
 for all continuous functions $\p:\,\R^d\to\R$ with a compact support.

The degree of concentration of $\g$ near diagonal of $\R^d\times\R^d$ shows the
closeness of $\nu_1$ and $\nu_2$ to  each other.

A continuous analogue of definition of uniformly spreading (\ref{d0}) (actually
 belonging to M.~Sodin, B.~Tsirelson \cite{ST}) has the form
 \begin{Def}
A locally finite positive measure $\nu$ on $\R^d$ is uniformly spread over $\R^d$, if
there is $\b>0$  such that
\begin{equation}\label{d1}
\Tra(\nu,\b\o)<\infty,
\end{equation}
where $\o$ is the Lebesgue measure on $\R^d$.
 \end{Def}
 In what follows we denote by $\nu^x$ the shift of the measure $\nu$ along $x\in\R^d$, and for any
$x=(x^1,\dots,x^d)\in\R^d$ put
 $$
 Q(x,r)=\{y=(y^1,\dots,y^d)\in\R^d:\,x^j-r/2\le y^j<x^j+r/2,\,j=1,\dots,d\}.
 $$
Also, we denote by $\chi^m(y),\,m\in\Z^d$, the indicator of the cube $Q(m,\,1)$.

For  a discrete sequence $(a_n)_{n\in\N}$  we set
\begin{equation}\label{a}
\nu=\sum_n\delta^{a_n},
\end{equation}
where $\delta$ is a unit mass sitting in the origin. Then
 \begin{equation}\label{i2}
C_1\Tra(\nu,\a^{-d}\o)\le \inf_\psi \sup_n |a_n-\a\psi(n)|\le C_2\Tra(\nu,\a^{-d}\o).
 \end{equation}
As above, the infimum  is taken over all bijections $\psi:\,\N\to\Z^d$, and the
constants $C_1$ and $C_2$ depend only on the dimension $d$.

In fact, if the sequence satisfies (\ref{d0}), then the measure
$\g=\sum\delta^{a_n}(x)\otimes\a^{-d}\chi^{\a\psi(n)}(y/\a)\o(y)$ is a transportation
measure between $\nu$ and $\a^{-d}\o$, and the first inequality in (\ref{i2}) follows
easily. The second inequality in (\ref{i2}) is nontrivial. Its proof in \cite{L2} is
based on the Rado Lemma from the graph theory.

  The main result of our article is the following theorem.
\begin{Th}\label{0}
A positive locally finite measure $\nu\not\equiv0$ is uniformly spread over $\R^d$ if
and only if there exists a constant $C_3<\infty$ such that
\begin{equation}\label{c}
  \Tra(\nu,\nu^z)<C_3\qquad\forall\,z\in\R^d.
\end{equation}
\end{Th}
{\bf Proof of the Theorem \ref{0}}.  By \cite[Theorem 1.2]{ST}, we have
  $$
\Tra(\nu_1,\nu_3)\le \Tra(\nu_1,\nu_2)+\Tra(\nu_2,\nu_3).
  $$
Hence, (\ref{d1}) and the equality $\o^z\equiv\o$ imply (\ref{c}).

Next, suppose that the measure $\nu$ satisfies (\ref{c}). We decompose $\R^d$ into the
cubes $Q(m,1),\,m\in\Z^d$. Let $\g_m$ be a transportation measure between $\nu$ and
$\nu^{m}$. For fixed $k\in\Z^d$ set
 $$
\l_{k,n}=\frac{1}{(2n+1)^d}\sum_{\|m-k\|_\infty\le n}\g_m.
 $$
It follows from definition of a transportation measure (\ref{t1}) and (\ref{t2}) that
$\l_{k,n}$ is a transportation measure between $\nu$ and
 $$
\mu_{k,n}=\frac{1}{(2n+1)^d}\sum_{\|m-k\|_\infty\le n}\nu^{m}.
 $$

Note that we can replace in (\ref{t1}) and (\ref{t2}) the function $\p$ by the
indicator function of any bounded Borel subset  of $\R^d$. Therefore,  by (\ref{c}),
we get for any $m,k\in\Z^d$
\begin{eqnarray*}
\nu^m(Q(k,1))&=&\nu^{m+k}(Q(0,1))=\g_{m+k}(Q(0,1)\times\R^d)\\&\le&\g_{m+k}(\R^d\times
Q(0,C_3+1))=\nu(Q(0,C_3+1)).
\end{eqnarray*}
Hence the measures  $\mu_{k,n}$ are uniformly bounded on every compact subset of
$\R^d$, and the measures $\l_{k,n}$ are uniformly bounded on every compact subset of
$\R^d\times\R^d$. By (\ref{c}),
\begin{equation}\label{t}
\Tra(\nu,\mu_{k,n})\le \sup_{m\in\Z^d}\Tra(\nu,\nu^{m})\le C_3.
\end{equation}
Take a subsequence $n'\to\infty$ such that for each $k\in\Z^d$ the measures $\l_{k,n}$
weakly converge to some measure $\l_{k}$, and the measures $\mu_{k,n}$ weakly converge
to some measure $\mu_{k}$. Note that for any $k,\,k'\in\Z^d$
 $$
 \mu_{k,n}-\mu_{k',n}=\frac{1}{(2n+1)^d}\left[\sum_{\|m-k\|_\infty\le n,\|m-k'\|_\infty>n}\nu^{m}-
\sum_{\|m-k'\|_\infty\le n,\|m-k\|_\infty>n}\nu^{m}\right],
 $$
 and
 $$
\card\{m:\,\|m-k\|_\infty\le n,\|m-k'\|_\infty>n\}=O(n^{d-1}),
 $$
 $$
 \card\{\|m-k'\|_\infty\le n,\|m-k\|_\infty>n\}=O(n^{d-1}),
 $$
as $n\to\infty$. Hence, for each $k,\,k'\in\Z^d$ the variations of measures
$\mu_{k,n}-\mu_{k',n}$ on every compact subset  tend to zero as $n\to\infty$ . The
same assertion is valid for the differences $\l_{k,n}-\l_{k',n}$. Therefore, for all
$k\in\Z^d$ we have  $\mu_k\equiv\mu$ and $\l_k\equiv\l$ for some measures $\l$ and
$\mu$. Moreover,
 $$
\mu(Q(k,1))=\mu(Q(k',1))\qquad\forall\, k,k'\in\Z^d.
 $$
It can easily be checked that  the measure $\l$ is a transportation measure between
$\nu$ and $\mu$. By (\ref{t}), we get
  $$
\Tra(\nu,\mu)\le C_3.
  $$
Hence, $\mu\not\equiv0$. Let $\rho_m$ be the restriction of the measure $\mu$ to
$Q(m,1)$. We obtain that the measure
 $$
\sum_{m\in\Z^d}\rho_m(x)\otimes\chi^{m}(y)\o(y)
 $$
is a transportation measure between $\mu$ and $\b\o$ with $\b=\mu(Q(0,1))$. Finally,
 $$
 \Tra(\nu,\b\o)\le \Tra(\nu,\mu)+\Tra(\mu,\b\o)\le C_3+1.
 $$
 Theorem \ref{0} is proved. \bs

\medskip
For discrete sequences in $\R^d$ we obtain the following result.
\begin{Th}
A discrete sequence $(a_n)_{n\in\N}\subset\Z^d$ satisfies (\ref{d0}) if and only if
for any $z\in\R^d$ there is a bijection $\sigma:\N\to\N$ such that
 \begin{equation}\label{s}
  \sup_n|a_n+z-a_{\sigma(n)}|\le C_7<\infty.
\end{equation}
\end{Th}

{\bf Proof}. Clearly, (\ref{d0}) yields (\ref{s}). On the other hand, if (\ref{s})
holds and $\nu$ is the measure defined in (\ref{a}), then the measure
 $$
\sum_{n\in\N}\delta^{a_n+z}(x)\otimes \delta^{a_{\sigma(n)}}(y)
 $$
is a transportation measure between $\nu^z$ and $\nu$. Condition (\ref{s}) implies
 $\Tra(\nu^z,\nu)<\infty$. By  Theorem \ref{0}, for some $\a>0$ we have
$\Tra(\nu,\a^{-d}\o)<\infty$. Using (\ref{i2}), we obtain (\ref{d0}). \bs

\bigskip

Mathematical School, Kharkov national university, Swobody sq.4, Kharkov, 61077 Ukraine

e-mail: Sergey.Ju.Favorov@univer.kharkov.ua

\medskip
Mathematics Division, Institute for Low Temperature Physics and Engineering, 47 Lenin
ave., Kharkov 61103, Ukraine

e-mail: artemdudko@rambler.ru

\end{document}